\documentclass[11pt]{amsart}
\usepackage{amscd,amssymb,amsmath}
\usepackage[arrow,matrix]{xy}
\usepackage{graphicx}
\usepackage{cite}
\usepackage{ifpdf}
\usepackage{tikz}

\newtheorem{thm}[subsection]{Theorem}
\newtheorem{lemma}[subsection]{Lemma}

\newtheorem{defn}[subsection]{Definition}
\newtheorem{ex}[subsection]{Example}

\numberwithin{equation}{section} 

\newcommand{\F}{{\mathcal F}}

\newcommand{\bea}{\begin{eqnarray}}
\newcommand{\eea}{\end{eqnarray}}

\begin{document}
\title [On evolution algebras]{On evolution algebras}

\author {J.M. Casas, M. Ladra, B.A. Omirov, U.A. Rozikov}

\address{J.\ M.\ Casas\\ Department of Applied Mathematics, E.U.I.T. Forestal, Pontevedra, University of Vigo, 36005, Spain.}
 \email {jmcasas@uvigo.es}
 \address{M.\ Ladra\\ Department of Algebra, University of Santiago de Compostela, 15782, Spain.}
 \email {manuel.ladra@usc.es}
 \address{B.\ A.\ Omirov and U.\ A.\ Rozikov\\ Institute of mathematics and information technologies,
Tashkent, 100125, Uzbekistan.} \email {omirovb@mail.ru,
rozikovu@yandex.ru}

\begin{abstract} The structural constants of an evolution
algebra is given by a quadratic matrix $A$. In this work we
establish equivalence between nil, right nilpotent evolution
algebras and evolution algebras, which are defined by upper
triangular matrix $A$. The classification of 2-dimensional complex
evolution algebras is obtained. For an evolution algebra with a
special form of the matrix $A$ we describe all its isomorphisms and
their compositions. We construct an algorithm running under Mathematica which  decides if two  finite dimensional evolution algebras are isomorphic.


\end{abstract}

\maketitle

\section{Introduction}
\large In this paper we consider a class of algebras called
evolution algebras. The concept of evolution algebra lies between
algebras and dynamical systems. Algebraically, evolution algebras
are non-associative Banach algebras; dynamically, they represent
discrete dynamical systems. Evolution algebras have many connections
with other mathematical fields including graph theory, group theory,
stochastic processes, mathematical physics, etc. \cite{Tian2},
\cite{Tian1}.

In the book \cite{Tian1} the foundation of evolution algebra theory
and applications in non-Mendelian genetics and Markov chains is
developed, with pointers to some further research topics.

Let $(E,\cdot)$ be an algebra over a field $K$. If it admits a basis
$e_1,e_2,\dots$, such that
$$e_i\cdot e_j=0, \ \ \mbox{if}\ \ i\ne j;$$
$$e_i\cdot e_i=\sum_{k}a_{ik}e_k, \ \ \mbox{for any}\ \ i,$$
then this algebra is called an {\it evolution algebra}. We denote by
$A=(a_{ij})$ the matrix of the structural constants of the evolution
algebra $E$.

In \cite{ly} an evolution algebra $\mathcal A$ associated to the
free population is introduced and using this non-associative algebra
many results are obtained in explicit form, e.g. the explicit
description of stationary quadratic operators, and the explicit
solutions of a nonlinear evolutionary equation in the absence of
selection, as well as general theorems on convergence to equilibrium
in the presence of selection.

In study of any class of algebras, it is important to describe up to
isomorphism even algebras of lower dimensions because such
description gives examples to establish or reject certain
conjectures. In this way in \cite{M} and \cite{Umlauf}, the
classifications of associative and nilpotent Lie algebras of low
dimensions were given.

In this paper we study some properties of evolution algebras. The
paper is organized as follows. In Section 2 we establish equivalence
between nil, right nilpotent evolution algebras and evolution
algebras, which are defined by upper triangular matrix $A$. In
\cite{Omirov} it was proved that these notions are equivalent to the
nilpotency of evolution algebras, but right nilindex and nilindex do
not coincide in general. Thus it is natural to study conditions when
some powers of the evolution algebras are equal to zero. In Section
3 we consider an evolution algebra $E$ with an upper triangular
matrix $A$ and drive a system of equation (for entries of the matrix
$A$) solutions to which gives $E^k=0$ for small values of $k$.
Section 4 is devoted to the classification of 2-dimensional complex
evolution algebras. In Section 5 for an evolution algebra with a
special form of the matrix $A$ we describe all its isomorphisms and
their compositions. Finally, in Appendix, we construct an algorithm running under Mathematica, using Gr\"obner bases and the star product of two evolution matrices, which  decides if two  finite dimensional evolution algebras are isomorphic.

\section{Nil and right nilpotent evolution algebras}

In this section we prove that notions of nil and right nilpotency
are equivalent for evolution algebras. Moreover, the defined matrix
$A$ of such algebras has upper (or lower, up to permutation of
basis of the algebra) triangular form.

\begin{defn} An element $a$ of an evolution algebra $E$  is called nil if
there exists $n(a) \in \mathbb{N}$ such that $(\cdots
((\underbrace{a\cdot a)\cdot a)\cdots a}_{n(a)})=0$.
Evolution algebra $E$ is called nil if every  element of the algebra
is nil.
\end{defn}

\begin{thm} \label{thm1} Let $E$ be a nil evolution algebra with basis $\{e_1, \dots, e_n\}$.
Then for the elements of the matrix $A=(a_{ij})$ the following
relation holds
\begin{equation}\label{1}
a_{i_1i_2}a_{i_2i_3}\dots a_{i_ki_1}=0 ,
\end{equation}
 for all
$i_1,\dots,i_k\in\{1,\dots, n\}$ and $k\in \{1,\dots,n\}$, with
$i_p\ne i_q$ for $p\ne q$.
\end{thm}
\begin{proof}
Note that $((e_i\cdot e_i)\cdot e_i)=a_{ii}e^2_i$, hence $a_{ii}=0$
(otherwise the element $e_i$ is not nil). We shall prove the
equality (\ref{1}) for right normed terms by induction.

For the element $e_i+e_j,  \ 1 \leq i, \ j \leq n$,  it can be
proved by induction the following relation
$$(e_i+e_j)^{2s}=a_{ij}^{s-1}a_{ji}^{s-1}(e_i+e_j)^2.$$
The nil condition for the element $e_i+e_j$ leads to the equality
for the elements of the matrix $A$: $$a_{ij}a_{ji}=0 \  \mbox{or}
\ e_i^2+e_j^2=0.$$

Take in account the fact that $a_{ii}=a_{jj}=0$ for any $i$ and $j$
and comparing the coefficients at the basic elements, from condition
$e_i^2+e_j^2=0$ we obtain $a_{ij}=a_{ji}=0$. Hence, the equation
$a_{ij}a_{ji}=0$ for all $i, \ j$ is obtained and therefore the
equality (\ref{1}) is true for $k=2$.

Let (\ref{1}) be true for $k-1$. We shall prove it for $k$. For this
purpose we consider element $e_{i_1}+e_{i_2}+\dots+e_{i_k}$. Without
loss of generality instead of this element we can consider the
following element $e_1+e_2+\dots+e_k$. Using the hypothesis of the
induction it is not difficult to note that
$$\left(\sum_{i=1}^k e_i\right)^{s+2}=
\sum_{\substack{i_1,\dots,i_s=1 \\ i_p\ne i_q, \ p\ne q}}^k a_{i_1
i_2}a_{i_2i_3}a_{i_3i_4}\dots a_{i_{s-1}i_s}e_{i_s}^2.$$

Let us take in the above expression $s=k+1$, then $i_{s-1}=i_k$.
From induction hypothesis the coefficient
$a_{i_1i_2}a_{i_2i_3}a_{i_3i_4}\dots a_{i_{s-1}i_s}$ is equal to
zero if $i_s\in\{i_2,\dots,i_{s-1}\}$. Therefore we need to consider
the case $i_s=i_1$ and the above expression will have the following
form
$$\left(\sum_{i=1}^k e_i\right)^{k+3}=\sum_{\phi\in S_k}
a_{\phi(1)\phi(2)}a_{\phi(2)\phi(3)}\dots
a_{\phi(k)\phi(1)}e^2_{\phi(1)}=$$
$$ \sum_{i=1}^k\left(\sum_{\phi\in S_k: \phi(1)=i}
a_{i\phi(2)}a_{\phi(2)\phi(3)}\dots a_{\phi(k)i}\right)e^2_i \, ,$$
where $S_k$ denotes the symmetric group of permutations of $k$ elements.

Denote
$$\F_i=\sum_{\phi\in S_k: \phi(1)=i}a_{i\phi(2)}a_{\phi(2)\phi(3)}\dots
a_{\phi(k)i} \, .$$ We need the following lemmas

\begin{lemma}\label{lem2} For any $i,j=1,\dots,k$ we have
$\F_i=\F_j$.
\end{lemma}
\proof For $\phi\in S_k$ with $\phi(1)=i$ we construct a unique
$\overline{\phi}\in S_k$ such that $\overline{\phi}(1)=j$ and
\begin{equation}\label{*}
a_{i\phi(2)}a_{\phi(2)\phi(3)}\dots
a_{\phi(k)i}=a_{j\overline{\phi}(2)}a_{\overline{\phi}(2)\overline{\phi}(3)}\dots
a_{\overline{\phi}(k)j}
\end{equation}
 as follows. Let $s$ be the number such that
$\phi(s)=j$, then the permutation $\overline{\phi}$ is defined as
$$\overline{\phi}=\left(\begin{array}{cccccccc}
1&2&\dots&k-s+1&k-s+2&k-s+3&\dots&k\\[2mm]
j&\phi(s+1)&\dots&\phi(k)&i&\phi(2)&\dots&\phi(s-1)\\[2mm]
\end{array}\right).$$
By construction, we note that for a given $\phi$ the
$\overline{\phi}$ is uniquely defined and satisfies (\ref{*}). Thus
we get $\F_i=\F_j$.
\endproof
Put $a=\sum_{i=1}^k e_i$.
\begin{lemma}\label{lem3} If $a^2=0$ then
$\F_1=0$.
\end{lemma}
\proof From $a^2=0$ we obtain
$$\sum_{\substack{i=1 \\ i\ne j}}^ka_{ij}=0, \ j=1,\dots,n.$$
Using this equality we get
$$\F_1=\sum_{\phi\in S_k: \phi(1)=1}a_{1\phi(2)}a_{\phi(2)\phi(3)}\dots
a_{\phi(k)1}=$$ $$ -\sum_{\phi\in S_k: \phi(1)=1}\sum_{\substack{i=2 \\ i\ne \phi(2)}}^k a_{i\phi(2)}a_{\phi(2)\phi(3)}\dots a_{\phi(k)1}.$$
Since for any $i=2,\dots,k$ there exists $s_i$ such that
$\phi(s_i)=i$, by the assumption of the induction we get
$$a_{i\phi(2)}a_{\phi(2)\phi(3)}\dots a_{\phi(k)1}=a_{i\phi(2)}a_{\phi(2)\phi(3)}\dots a_{\phi(s_i-1)i}a_{i\phi(s_i+1)}\dots
a_{\phi(k)1}=0.$$
\endproof

Now we continue the proof of theorem. Using Lemma \ref{lem2}, we get
$$\left(\sum_{i=1}^k e_i\right)^{k+3}=\F_1a^2=0.$$
By Lemma \ref{lem3} we get $\F_1=0$. Fix an arbitrary $\phi_0\in S_k$
with $\phi_0(1)=1$ and multiply both side of $\F_1=0$ by
$a_{1\phi_0(2)}a_{\phi_0(2)\phi_0(3)}\dots a_{\phi_0(k)1}$ then
(again using the assumption of the induction) we obtain
$$a^2_{1\phi_0(2)}a^2_{\phi_0(2)\phi_0(3)}\dots
a^2_{\phi_0(k)1}=0,$$ i.e.
$$a_{1\phi_0(2)}a_{\phi_0(2)\phi_0(3)}\dots
a_{\phi_0(k)1}=0,$$
 which completes
the induction and the proof of theorem.\end{proof}

For an evolution algebra $E$ we introduce the following sequence
$$E^{<1>}=E, \ \ \ E^{<k+1>}=E^{<k>}E, \ \  k \geq 1.$$
\begin{defn} An evolution algebra is called right nilpotent if there
exists some $s\in \mathbb{N}$ such that $E^{<s>}=0$.
\end{defn}

Let evolution algebra $E$ be a right nilpotent algebra, then it is
evident that $E$ is nil algebra. Therefore for the related matrix
$A=\left(a_{ij}\right)_{i,j=1}^n$ we have
\begin{equation}\label{11}
a_{i_1i_2}a_{i_2i_3}\dots a_{i_ki_1}=0,
\end{equation}
  for any $k\in \{1, 2, \dots,
n\}$ and arbitrary $i_1, i_2, \dots, i_k\in \{1,2, \dots, n\}$ with
$i_p\ne i_q$ for $p\ne q$.

\begin{lemma} \label{lem1} Let the matrix $A$ satisfies (\ref{11}). Then for any
$j\in \{1,\dots,n\}$ there is a row $\pi_j$ of $A$ with $j$ zeros.
Moreover, $\pi_{j_1}\ne \pi_{j_2}$ if $j_1\ne j_2$.
\end{lemma}
\begin{proof} First we shall prove that there is a row with
$j=n$ zeros, i.e. all zeros. Assume that there is not such a row.
Then for any $i\in\{1,\dots,n\}$ there is a number
$\beta(i)\in\{1,\dots,n\}\setminus \{i\}$ such that
$a_{i\beta(i)}\ne 0$. Consider the sequence
$$i_1=1, i_2=\beta(1),
i_3=\beta(\beta(1)),\dots,i_{n+1}=\underbrace{\beta(\cdots(\beta}_n(1))).$$
Then by assumption we have $a_{i_mi_{m+1}}\ne 0$, for all
$m=1,\dots,n$ hence
\begin{equation}\label{2}
a_{i_1i_2}a_{i_2i_3}\dots a_{i_ni_{n+1}}\ne 0.
\end{equation}
  Since $i_j\in
\{1,\dots,n\}, j=1,\dots,n+1$ between them there are $i_p=i_q$ for
some $p\ne q\in \{1,\dots,n+1\}$. Thus
\begin{equation}\label{3}
a_{i_pi_{p+1}}a_{i_{p+1}i_{p+2}}\dots a_{i_{q-1}i_p}\ne 0.
\end{equation}

 So (\ref{3}) is in contradiction with (\ref{11}). Thus there is a row $\pi_n$ with
 all zeros ($n$ zeros).

 Now we shall prove that there is a row $\pi_{n-1}\ne \pi_n$ of
 $A$  with $n-1$ zeros. Consider $A_{\pi_n}$-minor of $A$ which is
 constructed by $A$ deleting row $\pi_n$ and column $\pi_n$.
 Matrix $A_{\pi_n}$ is $(n-1)\times (n-1)$ and condition (\ref{11})
 implies the condition
 \begin{equation}\label{4}
a_{i_1i_2}a_{i_2i_3}\dots a_{i_ki_1}=0,
\end{equation}
 for any $k\in \{1,\dots,n-1\}$
and arbitrary $i_1,i_2,\dots,i_k\in \{1, \dots, n\}\setminus
\{\pi_n\}$ with $i_p\ne i_q$ for all $p\ne q$.

To prove that $A$ has a row with $j=n-1$ zeros it is enough to prove
that $A_{\pi_n}$ has a row with all zeros. But this problem is the
same as the case $j=n$, only we must consider condition (\ref{4})
instead of (\ref{11}). Iterating this argument we can show that for
any $j$ there exists a row $\pi_j$ with $j$ zeros. The proof of the
lemma is completed. \end{proof}

The following theorem is the main result of this section.
\begin{thm} \label{thm2} The following statements are equivalent for
an $n$-dimensional evolution algebra $E$:

a) The matrix corresponding to $E$ can be written as
\begin{equation}\label{5}
\widehat{A}=\left(
  \begin{array}{ccccc}
 0 & a_{12} & a_{13} &\dots &a_{1n} \\[1.5mm]
 0 & 0 & a_{23} &\dots &a_{2n} \\[1.5mm]
0 & 0 & 0 &\dots &a_{3n} \\[1.5mm]
\vdots & \vdots & \vdots &\cdots & \vdots \\[1.5mm]
0 & 0 & 0 &\cdots &0 \\
\end{array}
\right);
\end{equation}

b) $E$ is right nilpotent algebra;

c) $E$ is nil algebra.
\end{thm}
\begin{proof}

\textbf{ b) $\Rightarrow$ a).} Since the equality (\ref{11}) is true
for right nilpotent algebra then we are in the conditions of the
Lemma \ref{lem1}. Consider the permutation of the first indexes
$\{1, \dots, n\}$ of the matrix $A$ as
$$\pi=\left(
  \begin{array}{ccccc}
 1 & 2 & 3& \dots &n \\[1.5mm]
 \pi_1 & \pi_2 & \pi_3& \dots &\pi_n \\
\end{array}
\right),
$$
where $\pi_j$ is defined in the proof of the Lemma \ref{lem1}.
Note that Lemma \ref{lem1} is also true for columns: for any $j$
there is a column $\tau_j$ with $j$ zeros. Moreover $\tau_p\ne
\tau_q$, $p\ne q$. Now consider permutation of the second indexes
$\{1, \dots, n\}$ as
$$\tau=\left(
  \begin{array}{ccccc}
 1 & 2 & 3 & \dots &n \\[1.5mm]
 \tau_1 & \tau_2 & \tau_3& \dots &\tau_n \\
\end{array}
\right).
$$
Then $\tau(\pi(A))=\widehat{A}$.

The implication \textbf{ b) $\Rightarrow$ c)} is evident since every
right nilpotent evolution algebra is nil algebra.

The implication \textbf{ a) $\Rightarrow$ b, \ c )} is also true,
because the table of the multiplication of the evolution algebra
defined by upper triangular matrix $A$ will be right nilpotent and
nil.

The implication \textbf{ c) $\Rightarrow$ a)} follows from Theorem
\ref{thm1}. \end{proof}

\section{Conditions for $E^k=0$}

For an evolution algebra $E$ we define the \emph{lower central sequence by}:
$$E^1=E, \ \ \ E^k=\sum_{i=1}^{k-1}E^iE^{k-i}, \ k\geq 1.$$

An evolution algebra $E$ is called \emph{nilpotent} if there exists $n\in
\mathbb{N}$ such that $E^n=0$. In \cite{Omirov}, it is proved that
the notions of nilpotent and right nilpotent are equivalent.

In this section we consider an $n$-dimensional evolution algebra $E$
with a triangular (as $\widehat{A}$ in Theorem \ref{thm2}) matrix
$A$ and for small values of $k$ we present conditions on entries of
$A$ under which $E^k=0$.

First, for $n=3$ we have $E^2=0 \Leftrightarrow a_{ij}\equiv 0$
and $E^3=0 \Leftrightarrow a_{12}a_{23}=0$. For $n=4$ one easily
finds $E^2=0 \Leftrightarrow a_{ij}\equiv 0$ and
$$E^3=0 \Leftrightarrow
a_{12}a_{23}=0, a_{12}a_{24}=0, a_{13}a_{34}=0, a_{23}a_{34}=0.$$

Now we consider an arbitrary $n\in \mathbb{N}$ and for $k=3,4,5$, we
shall drive solutions of a system of equations (for $a_{ij}$) which
give $E^k=0, k=3,4,5$.

Let $E=<e_1, \dots, e_n>$ be an evolution algebra with matrix
(\ref{5}).

{\it Case $k=3$:}  We have $ e_i^2e_j=e_je_i^2, \ \forall i,j;
e_ie_je_k=0, i\ne j$. Thus we get
\[
e_i^2e_j=\begin{cases}
0 &  \mbox{if} \ \ j\leq i;\\
\sum^n_{s=j+1}a_{ij}a_{js}e_s &  \mbox{if} \ \ j\geq i+1, \,  s\geq
j+1 \, .\\
\end{cases}
\]
So the system of equations is
$$a_{ij}a_{js}=0\ \ \mbox{for} \ \ j\geq i+1, s\geq j+1, \,
i,j=1, \dots, n.\eqno(n;3)$$

{\it Case $k=4$:} Since $ e_i^2e_j^2=e_j^2e_i^2;\
(e_i^2e_j)e_s=(e_je_i^2)e_s=e_s(e_je_i^2)=e_s(e_i^2e_j); \
e_ie_je_se_t=0$, if $i\ne j$; $(e_ie_je_s)e_t=e_t(e_ie_je_s),\dots$
it will be enough to consider $e_i^2e_j^2$ and $(e^2_ie_j)e_s$. We
have
$$e_i^2e_j^2=\sum^n_{t=j+2}\left(\sum^{t-1}_{u=j+1}a_{iu}a_{ju}a_{ut}\right)e_t, \ \ i\leq j, \,
i,j=1,\dots,n;$$
$$(e_i^2e_j)e_s=\sum_{t=s+1}^na_{ij}a_{js}a_{st}e_t, \ \ j \ \geq i+1,\ s \geq j+1.$$
 So the system of equations is
$$\begin{cases}
\sum^{t-1}_{u=j+1}a_{iu}a_{ju}a_{ut}=0& \mbox{if} \ \ j\leq i, t\geq j+2, i,j=1, \dots, n;\\
a_{ij}a_{js}a_{st}=0& \mbox{if} \ \ j\geq i+1, s\geq j+1, t\geq s+1.
\end{cases}
\eqno(n;4)
$$

{\it Case $k=5$:} We should only use previous non-zero words and
multiply them to get a word of length 5:
$$e_i^2e_j^2e_s=\sum^n_{t=s+1}\left(\sum^{s-1}_{u=j+1}a_{iu}a_{ju}a_{us}a_{st}\right)e_t, \ \ i\leq
j,s\geq j+2,$$
$$(e_i^2e_j)e_s^2=a_{ij}\sum\limits_{u=s+2}^n \left(\sum_{t=s+1}^{u-1}a_{jt}a_{st}a_{tu} \right) e_u, \ \ j \leq i+1, \ s \geq j,$$
$$e_i^2e_je_se_v=\sum\limits_{u=v+1}^n\left(a_{ij}a_{js}a_{sv}a_{vu}\right)e_u, \ \ j
\geq i+1,\ s \geq j+1,\ v \geq s+1.$$
 Thus we get the following system of equations
$$\left\{\begin{array}{lll}
\sum^{s-1}_{u=j+1}a_{iu}a_{ju}a_{us}a_{st}=0\ \ \mbox{if} \ \ j\leq i, s\geq j+2, t\geq s+1,\\
a_{ij}\sum_{t=s+1}^{u-1}a_{jt}a_{st}a_{tu}=0\ \ \ \mbox{if} \ \
j\geq i+1, s\geq j, u\geq s+2,\\
a_{ij}a_{js}a_{sv}a_{vu}=0.
\end{array}
\right.\eqno(n;5)
$$
Thus we have proved the following
\begin{thm} Let $E$ be an evolution algebra with
matrix (\ref{5}) then $E^k=0$ if the elements of the matrix
(\ref{5}) satisfy the equations (n;k), where $k=3,4,5$.
\end{thm}

 \section{Classification of complex 2-dimensional evolution algebras}

In this section we give the classification of $2$-dimensional
complex evolution algebras.

Let $E$ and $E'$ be evolution algebras and $\{e_i\}$ a natural basis of $E$.
A linear map $ \varphi \colon E \to E'$ is called an \emph{homomorphism} of evolution algebras if it is an algebraic map and if the set $\{\varphi(e_i)\}$ can be complemented to a natural basis of $E'$.
Moreover, if $\varphi$ is bijective, then it is called an \emph{isomorphism}.

Let $E$ be a $2$-dimensional complex
evolution algebra and $\{e_1, e_2\}$ be a basis of the algebra
$E$.

It is evident that if $\dim E^2 =0$ then $E$ is an  abelian algebra,
i.e. an algebra with all products zero.

\begin{thm}\label{tt}
Any  2-dimensional complex evolution algebra $E$ is isomorphic to
one of the following pairwise non isomorphic algebras:
\begin{enumerate}
\item $\dim E^2=1$
\begin{itemize}
\item $E_1:\ \ e_{1}e_{1} = e_{1}$,

\item $E_2: \ \ e_{1}e_{1} = e_{1}, \ \ e_{2}e_{2}= e_{1}$,

\item $E_3: \ \ e_{1}e_{1} = e_{1} + e_{2}, \ \  e_{2}e_{2}= -e_{1}-
e_{2}$,

\item $E_4: \ \ e_{1}e_{1} = e_{2}$.
\end{itemize}
\item $\dim E^{2}=2$
\begin{itemize}

\item $E_5: \ \  e_{1}e_{1}=e_{1}+a_{2}e_{2}, \ \
e_{2}e_{2}=a_{3}e_{1}+e_{2}, \ \ 1 - a_{2}a_{3}\ne 0$, where
$E_5(a_{2},a_{3})\cong E_5'(a_{3},a_{2})$,

\item $E_6: \ \ e_{1}e_{1}=e_{2}, \ \ e_{2}e_{2}=e_{1}+a_{4}e_{2}$,
$a_4\ne 0$,  where $E_6(a_{4})\cong E_6(a'_{4}) \ \Leftrightarrow \
\frac{a'_4}{a_4}=\cos\frac{2\pi k}{3} + i \sin\frac{2\pi k}{3} \
\mbox{for some} \ k=0, 1, 2$.
\end{itemize}
\end{enumerate}
\end{thm}
\begin{proof} For an evolution algebra $E$ we have
$$e_{1}e_{1}= a_{1}e_{1}+a_{2}e_{2}, \ \  e_{2}e_{2}=
a_{3}e_{1}+a_{4}e_{2}, \ \ e_{1}e_{2}=e_{2}e_{1}=0.$$

Since $\dim E^{2} = 1$, then
$$e_{1}e_{1}=c_{1}(a_{1}e_{1}+a_{2}e_{2}), \ \
e_{2}e_{2}=c_{2}(a_{1}e_{1}+a_{2}e_{2}), \ \
e_{1}e_{2}=e_{2}e_{1}= 0.$$

Evidently  ($c_{1}, c_{2})\ne (0, 0)$, because otherwise our algebra
will be abelian.

Since $e_{1}$ and $e_{2}$ are symmetric we can suppose that
$c_{1}\ne 0$, then by simple change of basis (scale of it) we can
do $c_{1} =1$.

\textbf{Case 1.} $a_{1}\ne $ 0. Then we take the following change of
basis $$e'_{1} = a_{1}e_{1}+a_{2}e_{2 }, \ \  e'_{2} =
Ae_{1}+Be_{2},$$ where $a_{1}B-a_{2}A \ne 0$.

Consider the product
$$0 = e'_{1}e'_{2}
=(a_{1}e_{1}+a_{2}e_{2})(Ae_{1}+Be_{2})=a_{1}A(a_{1}e_{1}+a_{2}e_{2})
+$$ $$a_{2}Bc_{2}(a_{1}e_{1}+a_{2}e_{2})
=(a_{1}A+a_{2}Bc_{2})(a_{1}e_{1}+a_{2}e_{2}).$$ Therefore,
$a_{1}A+a_{2}Bc =0$, i.e. $A=-\frac{a_2 Bc_2 }{a_1 }$ and $
a_{1}B-a_{2}A=a_{1}B+\frac{a_2^2 Bc_2 }{a_1 }\ne 0$.

It means that in the case when $a^{2}_{1}+a^{2}_{2}c_{2 }\ne 0$ we
can take the above change.

Consider the products

$$e'_{1}e'_{1}=(a_{1}e_{1}+a_{2}e_{2})(a_{1}e_{1}+a_{2}e_{2})=a^{2}_{1}(a_{1}e_{1}+a_{2}e_{2})+
a^{2}_{2}c_{2}(a_{1}e_{1}+a_{2}e_{2})=$$
$$(a^{2}_{1}+a^{2}_{2}c_{2})(a_{1}e_{1}+a_{2}e_{2})=(a^{2}_{1}+a^{2}_{2}c_{2})e'_{1},$$
$$e'_{2}e'_{2}=(Ae_{1}+Be_{2})(Ae_{1}+Be_{2})=A^{2}(a_{1}e_{1}+a_{2}e_{2})+B^{2}_{2}c_{2}
(a_{1}e_{1}+a_{2}e_{2})=$$
$$(A^{2}+B^{2}c_{2})(a_{1}e_{1}+a_{2}e_{2})=\left(\frac{a_2^2 B^2c_2^2 }{a_1^2}+B^{2}c_{2}\right)e'_{1}=
\frac{B^2c_2 (a_1^2 +a_2^2 c_2 )}{a_1^2}e'_{1}.$$
\begin{description}
\item[\textbf{Case 1.1.}] $c_{2}=0$. Then $e_{1}e_{1}=a^{2}_{1}e_{1}, \ \
e_{2}e_{2}=e_{1}e_{2}=e_{2}e_{1}=0$. Taking $e'_{1}=\frac{e_1
}{a_1^2 }$ we get the algebra $E_1$.

\item[\textbf{Case 1.2.}] $c_{2 }\ne 0$. Then taking $B=\sqrt {\frac{a_1^2
}{c_2 }}$ we obtain
$$e_{1}e_{1}=(a^{2}_{1}+a^{2}_{2}c_{2})e_{1}, \ \ e_{2}e_{2}=(a^{2}_{1}+a^{2}_{2}c_{2})e_{1}.$$
If $a_1^2+a_2^2c_2\ne 0$, the following change of basis
$$e'_{1}=\frac{e_1 }{a_1^2 +a_2^2 c_2}, \ \ e'_{2}=\frac{e_2 }{a_1^2 +a_2^2c_2}$$
derives to the algebra with multiplication:
$$e_{1}e_{1}=e_{1}, \ \ e_{2}e_{2}=e_{1}.$$

If $a^{2}_{1}+a^{2}_{2}c_{2}=0$, then $c_2=-\frac{a_1^2}{a_2^2}$ and
we have $e_{1}e_{1}=a_{1}e_{1}+a_{2}e_{2}, \ \
e_{2}e_{2}=-\frac{a_1^3}{a_2^2}e_1-\frac{a_1^2}{a_2}e_2$.

The change of basis $e'_{1}=\frac{e_1}{a_1}, \ \
e'_{2}=\frac{a_2}{a_1^2}e_2$ derives to the algebra $E_3$.
\end{description}
\textbf{Case 2.} $a_{1}=0$. Then we have the products
$e_{1}e_{1}=a_{2}e_{2}, \ \ e_{2}e_{2}=c_{2}a_{2}e_{2}$, where
$a_{2}\ne 0$.

If $c_{2 }=0$, then by the change $e'_{1}=\frac{e_1}{\sqrt{a_2}}$ we
get again the algebra $E_4$.

If $c_{2 }\ne 0$, then by $e'_{1}=\frac{e_1}{\sqrt{c_2a^2_2}}, \ \
e'_{2}=\frac{e_2}{c_2a_2}$ we get the algebra
$e_{1}e_{1}=e_{2}, \ \ e_{2}e_{2}=e_{2}$ which is isomorphic to the
algebra $E_2$.

Now we consider algebras with $\dim E^{2}=2$. Let us write the table
of multiplication:

$$e_{1}e_{1}=a_{1}e_{1}+a_{2}e_{2}, \ \ e_{2}e_{2}=a_{3}e_{1}+a_{4}e_{2},$$
where $a_{1}a_{4 } - a_{2}a_{3}\ne 0$.

\textbf{Case 1.} $a_{1} \ne 0$ and $a_{4}\ne 0$. Then we can
transform both of them to unit, i.e. we can suppose $a_{1}=
a_{4}=1$. Therefore we have the two parametric family $E_5(a_{2},
a_{3})$:

$e_{1}e_{1}=e_{1}+a_{2}e_{2}, \ \ e_{2}e_{2}=a_{3}e_{1}+e_{2}, \ \
1 - a_{2}a_{3}\ne 0$.

Let us take the general change of basis of the form
$$e'_{1}=A_{1}e_{1}+A_{2}e_{2 }, \ \ e'_{2}=B_{1}e_{1} + B_{2}e_{2},$$ where $A_{1}B_{2} - A_{2}B_{1}\ne
0$.

Consider the product

$0=e'_{1}e'_{2}=(A_{1}e_{1}+A_{2}e_{2})(B_{1}e_{1}+B_{2}e_{2})=A_{1}B_{1}(e_{1}+a_{2}e_{2})
+ A_{2}B_{2}(a_{3}e_{1}+e_{2})=(A_{1}B_{1} +
A_{2}B_{2}a_{3})e_{1}+ (A_{1}B_{1}a_{2}+ A_{2}B_{2})e_{2}$.

Since in this new basis the algebra should be also evolution we have
$$A_{1}B_{1}+A_{2}B_{2}a_{3}=0, \ \ A_{1}B_{1}a_{2}+A_{2}B_{2}=0.$$
From which we have $A_{2}B_{2}(1 - a_{2}a_{3})=0, \ \ A_{1}B_{1}(1 -
a_{2}a_{3})=0$. Since $1 - a_{2}a_{3}\ne 0$, then we have
$A_{1}B_{1}=A_{2}B_{2}=0$.

\begin{description}
\item[\textbf{Case 1.1.}] $A_{2}=0$. Then $B_{1}=0$.

Consider the products

$e'_{1}e'_{1}=A^{2}_{1}(e_{1}+a_{2}e_{2})=e'_{1}+a'_{2}e'_{2}=A_{1}e_{1}+a'_{2}B_{2}e_{2}
\  \Rightarrow \  A^{2}_{1}=A_{1}, \  A^{2}_{1}a_{2}=a'_{2}B_{2} \
\Rightarrow  \ A_{1}=1$,

$e'_{2}e'_{2}=B^{2}_{2}(a_{3}e_{1}+e_{2})=a'_{3}e'_{1}+e'_{2}=a'_{3}A_{1}e_{1}+B_{2}e_{2}
\  \Rightarrow \  B^{2}_{2}a_{3}=a'_{3}A_{1}, \  B^{2}_{2}=B_{2} \
\Rightarrow \  B_{2}=1$.

\item[\textbf{Case 1.2.}] $A_{1}=0$. Then $B_{2}=0$ and from the family of
algebras $E_5(a_{2}, a_{3})$ we get the family $E_5(a_{3}, a_{2})$.
\end{description}
\textbf{Case 2.} $a_{1}=0$ or $a_{4}=0$. Since $e_{1}$ and $e_{2}$
are symmetric then without loss of generality we can suppose that
$a_{1}=0$.
$$e_{1}e_{1}=a_{2}e_{2}, \ \ e_{2}e_{2}=a_{3}e_{1}+a_{4}e_{2},$$ where $a_{2}a_{3}\ne
0$.

Taking the change of basis $e'_{1}=\sqrt[3]{\frac{1}{a_2^2a_3}}e_1,
\ \ e'_{2}=\sqrt[3]{\frac{1}{a_2a_3^2}}e_2$ we obtain the
one-parametric family of algebras $E_{6}(a_{4})$:
$$e_{1}e_{1}=e_{2}, \ \ e_{2}e_{2}=e_{1}+a_{4}e_{2}.$$

Let us take the general change of basis
$$e'_{1}=A_{1}e_{1}+A_{2}e_{2 }, \ \ e'_{2}=B_{1}e_{1}+B_{2}e_{2},$$ where $A_{1}B_{2} - A_{2}B_{1}\ne
0$.

Consider the product

$0=e'_{1}e'_{2}=(A_{1}e_{1}+A_{2}e_{2})(B_{1}e_{1}+B_{2}e_{2})=A_{1}B_{1}e_{2}+A_{2}B_{2}(e_{1}+a_{4}e_{2})$.

Therefore $A_{1}B_{1}+A_{2}B_{2}a_{4}=0, \ \ A_{2}B_{2}=0 \
\Rightarrow \ A_{1}B_{1}=0, \ A_{2}B_{2}=0$.

Without loss of generality we can assume that $A_{2}=0$. Then
$B_{1}=0$.

Consider the product

$e'_{1}e'_{1}=A^{2}_{1}e_{2}=e'_{2}=B_{2}e_{2} \ \Rightarrow \
 A^{2}_{1}=B_{2}$,

$e'_{2}e'_{2}=B^{2}_{2}(e_{1}+a_{4}e_{2})=e'_{1}+a'_{4}e'_{2}=A_{1}e_{1}+a'_{4}B_{2}e_{2}
\ \Rightarrow \ B^{2}_{2}=A_{1}, \ \  B^{2}_{2}a_{4}=B_{2}a'_{4}$.

From these equalities we have $B^{3}_{2}=1, \ \ B_{2}a_{4}=a'_{4}$.

If $\frac{a'_4}{a_4}=\cos\frac{2\pi k}{3} + i \sin\frac{2\pi k}{3} \
\mbox{for some} \ k=0, 1, 2$, then putting $B_{2}=\cos\frac{2\pi
k}{3} + i \sin\frac{2\pi k}{3}$ we obtain the isomorphism between
algebras $E_6(a_{4})$ and $E_6(a'_{4})$.

The pairwise non isomorphic
obtained algebras can be checked by comparison of the algebraic properties listed in the following table.
\begin{center}
\begin{tabular}{|l|c|c|c|c|}
\hline   &  $\dim E^2$ & Right Nilpotency  & $\dim$(Annihilator) & Nil Elements  \\
\hline $E_1$ & 1 & No & 1 &  Yes\\
 \hline $E_2$ & 1 & No & 0 & Yes\\
\hline $E_3$ & 1 & No  & 0 & No\\
\hline $E_4$ & 1 & Yes &  1&  Yes   \\
\hline $E_5$ & 2 & No & 0 & Non \\
\hline $E_6$ & 2 & No & 0 & Yes \\
\hline
\end{tabular}
\end{center}
\end{proof}

 \section{Isomorphisms of evolution algebras}

Since the study of the isomorphisms for any class of algebras is a crucial task and
taking into account the great difficulties of their description, in this section we consider a particular case of evolution algebras, which have
matrices in the diagonal $2\times 2$
non-zero blocks.

Let $E$ be an evolution algebra which has a matrix $A$ in the following form
$$A=\left(
  \begin{array}{ccccccc}
 a_1 & b_1 & 0 & 0 &\dots & 0 & 0 \\[1.5mm]
 c_1 & d_1 & 0 & 0 &\dots & 0 & 0 \\[1.5mm]
 0 & 0 & a_2 & b_2 &\dots & 0 & 0 \\[1.5mm]
 0 & 0 & c_2 & d_2 &\dots & 0 & 0 \\[1.5mm]
 \vdots &\vdots & \vdots & \vdots &\dots &\vdots&\vdots \\[1.5mm]
0 & 0 & 0 & 0 &\dots & a_n & b_n \\[1.5mm]
 0 & 0 & 0 & 0 &\dots &c_n & d_n \\[1.5mm]
\end{array}
\right). $$ The basis $\{e_1, \dots, e_{2n}\}$ of this evolution
algebra has the following relations:
$$e_ie_j=0, i\ne j; \ \ e_{2k-1}^2=a_ke_{2k-1}+b_ke_{2k}, \ \ k=1,2,\dots,n; $$
$$e_{2k}^2=c_ke_{2k-1}+d_ke_{2k}, \ \ k=1,2,\dots,n. $$
Let $\varphi$ be an isomorphism of the evolution algebra $E$ onto $E$ with matrix $A'$. Write $\varphi$ as
$$\varphi=\left(
  \begin{array}{cccc}
 \alpha_{11} & \alpha_{12} &\dots & \alpha_{1,2n} \\[1.5mm]
 \alpha_{21} & \alpha_{22} &\dots & \alpha_{2,2n} \\[1.5mm]
 \vdots & \vdots & \dots &\vdots \\[1.5mm]
\alpha_{2n,1} & \alpha_{2n,2} &\dots & \alpha_{2n,2n} \\[1.5mm]
\end{array}
\right), $$ with $\det(\varphi)\ne 0$. We have
$$(e'_i)^2=(\varphi(e_i))^2=(\alpha_{i1}^2a_1+\alpha_{i2}^2c_1)e_1+
(\alpha_{i1}^2b_1+\alpha_{i2}^2d_1)e_2+$$$$+
\dots+(\alpha_{i,2n-1}^2a_n+\alpha_{i,2n}^2c_n)e_{2n-1}+
(\alpha_{i,2n-1}^2b_n+\alpha_{i,2n}^2d_n)e_{2n}; \ \ i=1, 2, \dots,
2n.$$  For $i\ne j$ we get
\begin{equation}\label{10}\begin{array}{ll}
e_i'e_j'=\left(\alpha_{i1}\alpha_{j1}a_1+\alpha_{i2}\alpha_{j2}c_1\right)e_1+
\left(\alpha_{i1}\alpha_{j1}b_1+\alpha_{i2}\alpha_{j2}d_1\right)e_2+ \cdots\\
 +(\alpha_{i,2n-1}\alpha_{j,2n-1}a_n+\alpha_{i,2n}\alpha_{j,2n}c_n )e_{2n-1} \left(\alpha_{i,2n-1}\alpha_{j,2n-1}b_n+\alpha_{i,2n}\alpha_{j,2n}d_n\right)e_{2n}=0.\end{array}
\end{equation}
From (\ref{10}) we obtain
\begin{equation}\label{111}
\left\{\begin{array}{lllll}
\alpha_{i1}\alpha_{j1}a_1+\alpha_{i2}\alpha_{j2}c_1=0\\
\alpha_{i1}\alpha_{j1}b_1+\alpha_{i2}\alpha_{j2}d_1=0\\
\qquad \dots \\
\alpha_{i,2n-1}\alpha_{j,2n-1}a_n+\alpha_{i,2n}\alpha_{j,2n}c_n=0\\
\alpha_{i,2n-1}\alpha_{j,2n-1}b_n+\alpha_{i,2n}\alpha_{j,2n}d_n=0.\\
\end{array}\right.
\end{equation}
Let $S_{2n}$ be the group of permutations of ${1, 2, \dots, 2n}$.
\begin{thm} Assume that $\det(A)\ne 0$ then

(i) For any isomorphism $\varphi:E\to E$ there exists unique
$\pi=\pi(\varphi)\in S_{2n}$ such that
$$\varphi\in \Phi_\pi=\left\{\left(
  \begin{array}{cccc}
 \alpha_{11} & \dots & \alpha_{1,2n} \\[1.5mm]
 \alpha_{21} & \dots & \alpha_{2,2n} \\[1.5mm]
 \vdots &  \dots &\vdots \\[1.5mm]
\alpha_{2n,1} &\dots & \alpha_{2n,2n} \\[1.5mm]
\end{array}
\right): \begin{array}{cc} \alpha_{i\pi(i)}\ne 0, \ 1 \leq i
\leq2n & \\[1.5mm] \mbox{and the rest of elements} \
\alpha_{ij}=0 &\\[1.5mm]\end{array}\right\}.$$ Moreover, the set $\Phi=\cup_{\pi\in
S_{2n}}\Phi_\pi$ is the set of all possible homomorphisms.

(ii) For any $\pi,\tau\in S_{2n}$ the following equality holds
$$\Phi_\pi\Phi_\tau=\{\varphi\psi: \varphi\in \Phi_\pi, \psi\in \Phi_\tau\}=\Phi_{\tau\pi}.$$
The set $G=\{\Phi_\pi: \pi\in S_{2n}\}$ is a multiplicative group.
\end{thm}
\begin{proof}
(i) Since $\det(A)\ne 0$ we have $a_id_i-b_ic_i\ne 0$ for any
$i=1,2,\dots,n$. Thus from (\ref{111}) we have
\begin{equation}\label{12}
\alpha_{ik}\alpha_{jk}=0, \ i\ne j, \ i,j,k=1,\dots,2n.
\end{equation}
 By (\ref{12}) it is
easy to see that each row and each column of the matrix $\varphi$
must contain exactly one non-zero element. It is not difficult to
see that every such matrix $\varphi$ corresponds to a permutation
$\pi$. The set of all possible solutions of (\ref{12}) give all the
possible isomorphisms, i.e. we get the set $\Phi$.

(ii) Take $\varphi=\{\alpha_{ij}\}\in \Phi_\pi$ and
$\psi=\{\beta_{ij}\}\in \Phi_\tau$. Denote $\varphi \circ
\psi=\{\gamma_{ij}\}$. It is easy to see that
$$\gamma_{ij}=\begin{cases}
0& \mbox{if} \ \ j\ne \tau(\pi(i));\\
\alpha_{i\pi(i)}\beta_{\pi(i)\tau(\pi(i))}& \mbox{if} \ \
j=\tau(\pi(i)).\\
\end{cases}$$
This gives $\Phi_\pi\Phi_\tau=\Phi_{\tau\pi}$
and then one easily can check that $G$ is a group.
\end{proof}
Now for a fixed $\varphi$ (i.e. $\pi$) we shall find the matrix
$A'$. Consider $\pi\in S_{2n}$ and the corresponding
$\varphi_{\pi}=(\alpha_{ij})$:
$$\alpha_{ij}=\begin{cases}
0& \mbox{if} \ \ j\ne \pi(i);\\
\alpha_{i\pi(i)}& \mbox{if} \ \
j=\pi(i).\\
\end{cases}$$
We have
\begin{equation}\label{13}
e'_i=\alpha_{i\pi(i)}e_{\pi(i)}, \ \ i=1,\dots,2n.
\end{equation}
 Using this
equality we get
\[(e'_i)^2=\alpha^2_{i\pi(i)}e^2_{\pi(i)}=\begin{cases}
\alpha^2_{i(2k-1)}(a_ke_{2k-1}+b_ke_{2k}), & \ \ \mbox{if} \ \
\pi(i)=2k-1 \, ;\\[2mm]
\alpha^2_{i(2k)}(c_ke_{2k-1}+d_ke_{2k}),  &  \ \ \mbox{if} \ \
\pi(i)=2k \, .\\
\end{cases}
\]
By (\ref{13}) from the last equality we get
$$(e'_i)^2=\left\{\begin{array}{ll}
(\alpha_{i(2k-1)}a_k)e'_i+(\frac{\alpha^2_{i(2k-1)}}
{\alpha_{i(2k)}}b_k)e'_{\pi^{-1}(2k)} \, , \ \ \mbox{if} \ \
\pi(i)=2k-1 \, ; \\[2mm]
(\frac{\alpha^2_{i(2k)}}
{\alpha_{i(2k-1)}}c_k)e'_{\pi^{-1}(2k-1)}+(\alpha_{i(2k)}d_k)e'_{i} \, ,
\ \ \mbox{if} \ \
\pi(i)=2k.\\
\end{array}\right.$$
Thus $A'=(a'_{ij})$ is a matrix with
\begin{equation}\label{15}
a'_{ij}=\begin{cases} \alpha_{i(2k-1)}a_k  & \ \
\mbox{if} \ \
\pi(i)=2k-1, \ j=i \, ; \\[2mm]
\frac{\alpha^2_{i(2k-1)}} {\alpha_{i(2k)}}b_k,  & \ \ \mbox{if} \ \
\pi(i)=2k-1, \ \pi(j)=2k \, ; \\[2mm]
\frac{\alpha^2_{i(2k)}}{\alpha_{i(2k-1)}}c_k   & \ \ \mbox{if} \ \
\pi(i)=2k, \ \pi(j)=2k-1 \, ;\\[2mm]
\alpha_{i(2k)}d_k,   & \ \ \mbox{if} \ \
\pi(i)=2k, \ j=i \, ; \\[2mm]
0  &  \ \ \mbox{otherwise.}
\end{cases}
\end{equation}

\begin{thm} Assume that $\det(A)\ne 0$. Let $\varphi : E \to E$ ($A\to A'$)
be an isomorphism  then $A'$ has the same form as $A$ if and only if
$\varphi$ belongs to $\Phi_\pi$, where $\pi\in
S^b_{2n}=\{\pi=(\pi(1),\dots,\pi(2n))\in S_{2n}: \pi(i)\in
\{\pi(i-1)\pm 1\}, \ i=1, 2, \dots, 2n\}$.
\end{thm}
\begin{proof} Using given above formula (\ref{15}) for $A'$ and the condition $\det(A)\ne 0$ one can see that
it has form as $A$ iff $\pi(i)\in \{\pi(i-1)\pm 1\}, \ i=1, 2,
\dots, 2n\ $.
\end{proof}

Properties of the matrix $A$ can be uniquely defined by properties
of its non-zero blocks. So if we consider $n=1$ then for $\det(A)\ne
0$  we have two classes of isomorphisms:
$$ \Phi_{12}=\left\{\left(
  \begin{array}{cc}
 \alpha & 0 \\[1.5mm]
 0 & \delta \\[1.5mm]
 \end{array}
\right): \alpha\delta\ne 0\right\} \ \ \mbox{with} \ \ A'=\left(
  \begin{array}{cc}
 a\alpha & b{\alpha^2\over \delta}  \\[1.5mm]
c{\delta^2\over \alpha} & d\delta \\[1.5mm]
 \end{array}
\right).$$
$$ \Phi_{21}=\left\{\left(
  \begin{array}{cc}
 0 & \beta \\[1.5mm]
 \gamma & 0 \\[1.5mm]
 \end{array}
\right): \beta\gamma\ne 0\right\} \ \ \mbox{with} \ \ A'=\left(
  \begin{array}{cc}
 d\beta & c{\beta^2\over \gamma} \\[1.5mm]
 b{\gamma^2\over \beta} & a\gamma \\[1.5mm]
 \end{array}
\right).$$

It is easy to check the following embedding:
$$ \Phi_{12}\Phi_{21}\subset \Phi_{21}, \ \ \Phi_{21}\Phi_{12}\subset
\Phi_{21};$$
$$ \Phi_{21}\Phi_{21}\subset \Phi_{12}, \ \ \Phi_{12} \ \ \mbox{is a
group}.$$

Adding the symmetric property to the matrices $A$ and $A'$ we get
the following classes of isomorphisms:
$$ \Phi_{12}^s=\left\{\left(
  \begin{array}{cc}
 \alpha & 0 \\[1.5mm]
 0 & \delta \\[1.5mm]
 \end{array}
\right): \alpha\delta\ne 0, ad-b^2\ne 0,
b\alpha^3=b\delta^3\right\},$$
$$ \Phi_{21}^s=\left\{\left(
  \begin{array}{cc}
 0 & \beta \\[1.5mm]
 \gamma & 0 \\[1.5mm]
 \end{array}
\right): \beta\gamma\ne 0, ad-b^2\ne 0,
b\beta^3=b\gamma^3\right\}.$$

\section*{Appendix}
The following program written in Mathematica permits to check the existence (or non existence) of an isomorphism between two
evolution  algebras of dimension $n$. It is based on the star product of two evolution matrices, $A*B$, see \cite[page 31]{Tian2} and the computation of Gr\"obner bases. In particular, one can check again that the algebras $E_i, i=1,\dots,6$ (see Theorem \ref{tt}) are pairwise non isomorphic.

\

  \begin{verbatim}
StarProduct[A_List, B_List] := Module[{icont, jcont, kcont, 
AEB, ndim, Indices}, ndim = Dimensions[A][[1]];
   Indices = {};
   Do[Indices = Join[Indices, {{icont, jcont}}];
     ,{icont, 1, ndim}, {jcont, icont + 1, ndim}];
    AEB = Table[Table[aux[icont, jcont], {icont, 1, ndim}],
                {jcont,1, (ndim^2 - ndim)/2}];
     Do[AEB[[icont, kcont]] =
     A[[Indices[[icont]][[1]], kcont]]*
      B[[Indices[[icont]][[2]], kcont]];
       ,{icont, 1, Length[Indices]}, {kcont, 1, ndim}];
        Return[AEB];]
   SystemEquations[P_List, Q_List] :=Module[{FirstEquation,
   SecondEquation, ThirdEquation, A, ndim, Result},
    ndim = Dimensions[P][[1]];
    A = Table[Table[aux[icont, jcont], {jcont, 1, ndim}],
     {icont, 1, ndim}];
    FirstEquation = (A*A).Q - P.A;
    SecondEquation = StarProduct[Transpose[A], Transpose[A]].Q;
    ThirdEquation = {Det[A]*Y - 1};
    Result = Join[Flatten[FirstEquation], 
    Flatten[SecondEquation], ThirdEquation]; Return[Result];]
    IsoEvolAlgebrasQ[A1_, A2_] := Module[{Equations, BGrobner},
    Equations = SystemEquations[A1, A2];
    BGrobner = GroebnerBasis[Equations, Variables[Equations]];
    (* Print temporal *)
    Print[BGrobner];
    If[BGrobner == {1},
    Print["Evolution algebras are NOT isomorphic"];
      Print["Evolution algebras are isomorphic"]; ];]
\end{verbatim}

\

\begin{ex}
We check that the evolution algebras $E_5$ and $E_6$ are not isomorphic.
\begin{verbatim}
IsoEvolAlgebrasQ[{{1, a2}, {a3, 1}}, {{0, 1}, {1, a4}}]
{1}
Evolution algebras are NOT isomorphic
\end{verbatim}
\end{ex}

\section*{ Acknowledgements}

 The first and second authors were supported by Ministerio
de Ciencia e Innovaci\'on (European FEDER support included), grant
MTM2009-14464-C02, and by Xunta de Galicia, grant Incite09 207 215
PR. The third author was partially supported by the Grant
NATO-Reintegration ref. CBP.EAP.RIG. 983169. The fourth author
thanks to the Department of Algebra, University of Santiago de
Compostela, Spain,  for providing financial support of his visit to
the Department.

\end{document}